\documentclass[11pt]{article}
\usepackage[T1]{fontenc}
\usepackage[utf8]{inputenc}
\usepackage{lmodern,microtype}
\usepackage{amsmath,amssymb,amsthm,mathtools,bm}
\usepackage[margin=1in]{geometry}
\usepackage{booktabs,tabularx,array,enumitem}
\usepackage[hidelinks]{hyperref}
\usepackage[nameinlink,capitalise,noabbrev]{cleveref}
\hypersetup{
  pdftitle={Cutoff and a Gaussian-shift profile for the Rudvalis shuffle},
  pdfauthor={Yunjiang Jiang},
  pdfsubject={Rudvalis shuffle, mixing time, cutoff, Gaussian profile, interchange process}
}
\setlist[itemize]{topsep=3pt,itemsep=2pt,parsep=0pt}
\newtheorem{theorem}{Theorem}[section]
\newtheorem{lemma}[theorem]{Lemma}
\newtheorem{proposition}[theorem]{Proposition}
\newtheorem{corollary}[theorem]{Corollary}

\theoremstyle{remark}
\newtheorem{remark}[theorem]{Remark}

\newcommand{\E}{\mathbb E}
\newcommand{\Pp}{\mathbb P}
\newcommand{\C}{\mathbb C}
\newcommand{\R}{\mathbb R}
\newcommand{\Sn}{\mathfrak S_n}
\newcommand{\Tr}{\operatorname{Tr}}
\newcommand{\TV}{\operatorname{TV}}
\newcommand{\HS}{\operatorname{HS}}

\newcommand{\erf}{\operatorname{erf}}
\newcommand{\id}{\mathrm{id}}
\newcommand{\ind}{\mathbf 1}
\newcommand{\norm}[1]{\left\lVert #1\right\rVert}
\newcommand{\ip}[2]{\left\langle #1,#2\right\rangle}
\newcommand{\floor}[1]{\left\lfloor #1\right\rfloor}
\newcommand{\ceil}[1]{\left\lceil #1\right\rceil}
\newcommand{\ee}{\mathrm e}

\title{Cutoff and a Gaussian-shift profile for the Rudvalis shuffle}
\author{Yunjiang Jiang}
\date{October 3, 2026}

\begin{document}
\maketitle

\begin{abstract}
The Rudvalis shuffle removes the top card and reinserts it uniformly at random in one of the bottom two positions.  Hildebrand proved the classical upper bound $O(n^3\log n)$, and Wilson later proved the sharp lower coefficient $1/(8\pi^2)$, establishing the correct order but leaving the leading upper constant, cutoff, and the cutoff profile open.  We prove that the full labelled shuffle has total-variation cutoff at
\[
  \frac{n^3\log n}{8\pi^2}
\]
with a window of order $n^3$.  More precisely, for every fixed $s\in\R$,
\[
 d_n\!\left(\floor{\frac{n^3}{8\pi^2}(\log n+2s)}\right)
 \longrightarrow
 \erf\!\left(\frac{\ee^{-s}}2\right),
\]
and at the same times the chi-square divergence converges to
$\exp(2\ee^{-2s})-1$.

The upper bound is obtained by writing one deterministic rotation of the shuffle as a nonnormal sweep of adjacent averaging projections, comparing that sweep to a weighted interchange generator, and controlling the entire regular-representation heat trace.  The key sharp-window ingredients are heat-kernel smoothing through the Octopus Inequality, zero extension of deletion-harmonic ordered-tuple functions, a representation-theoretic multiplicity payment, and a sectorwise entrance estimate that disposes of all large Young-diagram deficits.  For the profile, an exact complex one-card eigenfunction yields a deck statistic whose surviving first Fourier mode converges to a Gaussian shift.  A matching full-deck $L^2$ bound then shows that this one-dimensional Gaussian likelihood ratio asymptotically exhausts all remaining information, giving the exact total-variation profile.
\end{abstract}



\tableofcontents
\newpage

\section{Introduction and historical context}\label{sec:intro}

\subsection{The shuffle and the long-standing sharp-constant problem}
The Rudvalis shuffle is one of the simplest examples of a card shuffle whose global mixing is much slower than the visible motion of a single step might suggest.  From a deck of $n$ labelled cards, remove the top card and place it either at the bottom or immediately above the bottom card, each with probability $1/2$.  The model was proposed by Arunas Rudvalis and was recorded in Diaconis's monograph on group representations and probability~\cite{Diaconis1988}.  We number positions from top to bottom.

The early quantitative theory was developed by Hildebrand, whose thesis established an $O(n^3\log n)$ total-variation upper bound~\cite{Hildebrand1990}.  Diaconis and Saloff-Coste subsequently developed comparison methods for random walks on finite groups and studied closely related symmetric variants of the Rudvalis dynamics~\cite{DSC1993,DSC1995}.  Wilson's 2003 paper supplied the decisive lower-bound advance: for every fixed $\varepsilon>0$, at time
\[
  \frac{1-o(1)}{8\pi^2}n^3\log n
\]
the original Rudvalis shuffle remains at total-variation distance at least $1-\varepsilon$ from uniform~\cite[Theorem~1]{Wilson2003}.  Thus the order $\Theta(n^3\log n)$ was known, but the upper and lower leading constants did not match.  Wilson's method also clarified why the problem is delicate: the relevant one-card mode is complex and its phase rotates quickly, while its modulus decays only on the much longer scale $n^3$~\cite{Wilson2003}.  Jonasson later gave a streamlined route to the same order of lower bound using an $n$-step viewpoint~\cite{Jonasson2006}, and Goel studied the broader family of top-to-bottom-$k$ shuffles~\cite{Goel2006}.

A different line of work studies hydrodynamic and fluctuation limits of particle projections of Rudvalis-type shuffles.  Gon\c{c}alves, Jara, Marinho and Moreira obtained such scaling limits and, in their 2023 paper, explicitly noted that matching mixing-time bounds for the original Rudvalis shuffle were still unknown~\cite{GoncalvesEtAl2023}.  Their work also emphasizes a useful conceptual separation: diffusive or hydrodynamic information about a fixed projection does not by itself control the full labelled permutation.

The present argument closes the sharp upper-bound gap and goes further to identify the cutoff window and profile.  A recent ingredient is the use of the Octopus Inequality of Caputo, Liggett and Richthammer~\cite{CLR}, in a heat-kernel smoothing comparison of the type that also appears in recent work of Hermon, Kozma and Salez on the exclusion process~\cite{HKS}.  We use only the comparison inequality derived from the Octopus Inequality; no exclusion-process cutoff theorem is imported into the Rudvalis shuffle.

\subsection{Main result}
Let $\pi_n$ be the uniform probability measure on $\Sn$, let $\mu_t$ be the law of the deck after $t$ elementary Rudvalis moves from a deterministic initial deck, and define
\[
 d_n(t)=\norm{\mu_t-\pi_n}_{\TV}.
\]
Because the shuffle is a random walk on the group $\Sn$, translation of the initial deck merely translates the law, and $d_n(t)$ does not depend on the deterministic starting permutation.

For probability measures $\mu\ll\pi_n$, define the chi-square divergence
\[
 \chi^2(\mu,\pi_n)
 =\E_{\pi_n}\!\left[\left(\frac{d\mu}{d\pi_n}-1\right)^2\right].
\]
We abbreviate $\chi_n^2(t)=\chi^2(\mu_t,\pi_n)$.

\begin{theorem}[Cutoff profile]\label{thm:mainprofile}
For each fixed $s\in\R$, put
\[
 t_n(s)=\floor{\frac{n^3}{8\pi^2}(\log n+2s)}.
\]
Then
\begin{align}
 d_n(t_n(s))&\longrightarrow
  2\Phi\!\left(\frac{\ee^{-s}}{\sqrt2}\right)-1
  =\erf\!\left(\frac{\ee^{-s}}2\right),\label{eq:main-tv}\\
 \chi_n^2(t_n(s))&\longrightarrow
  \exp(2\ee^{-2s})-1.\label{eq:main-chi}
\end{align}
Here $\Phi$ is the standard normal distribution function and
$\erf(x)=\frac2{\sqrt\pi}\int_0^x e^{-u^2}\,du$.
\end{theorem}

The profile is best understood as a \emph{Gaussian-shift profile}.  The total-variation profile itself is not a Gaussian density.  At the cutoff window, one real Fourier statistic has limiting law $N(0,1)$ under stationarity and $N(\sqrt2\ee^{-s},1)$ under the shuffled deck.  The total-variation distance between those two normal laws is exactly the right side of~\eqref{eq:main-tv}.

\begin{corollary}[Sharp mixing time and window]\label{cor:mixing}
For every fixed $\varepsilon\in(0,1)$,
\[
 t_{\mathrm{mix}}(\varepsilon)
 =\frac{n^3\log n}{8\pi^2}
 -\frac{n^3}{4\pi^2}
   \log\!\left(2\erf^{-1}(\varepsilon)\right)
 +o(n^3).
\]
In particular,
\[
 t_{\mathrm{mix}}(\varepsilon)
 =\left(\frac1{8\pi^2}+o(1)\right)n^3\log n,
\]
and the cutoff window has order $n^3$.  With the conventional centering
$t_n=n^3\log n/(8\pi^2)$, a natural window scale is
$w_n=n^3/(4\pi^2)$.
\end{corollary}

\begin{proof}
The function $F(s)=\erf(\ee^{-s}/2)$ is continuous and strictly decreasing, with $F(s)\to1$ as $s\to-\infty$ and $F(s)\to0$ as $s\to+\infty$.  Solving $F(s)=\varepsilon$ gives
$s=-\log(2\erf^{-1}(\varepsilon))$.  Monotonicity of total-variation distance in time and \cref{thm:mainprofile} then bracket the level-$\varepsilon$ hitting time between any two nearby fixed window parameters, yielding the stated expansion.
\end{proof}

\subsection{Why the constant $1/(8\pi^2)$ appears}
The proof makes the constant transparent.  One full turn of the deterministic rotation consists of $n$ elementary moves and becomes, in a rotating coordinate frame, a product
\[
 K=P_{n-1}\cdots P_0,
 \qquad P_i=\frac{I+U_{(i,i+1)}}2,
\]
of adjacent averaging projections around the cycle.  We compare this nonnormal product to the positive weighted interchange operator
\[
 L_Q=Q^{-1}A_0+\sum_{i=1}^{n-1}A_i,
 \qquad A_i=I-P_i.
\]
The endpoint comparison has the asymptotically sharp constant
\[
 K^*L_QK\preceq c_Q(I-K^*K),
 \qquad c_Q=\frac{1+5/Q}{4}\longrightarrow\frac14.
\]
On the other hand the first nonzero one-card eigenvalue of $L_Q$ is asymptotic to $2\pi^2/n^2$.  Hence the critical \emph{static} heat-trace time is $n^2\log n/(2\pi^2)$.  Multiplication by the endpoint factor $1/4$ gives $n^2\log n/(8\pi^2)$ sweeps, and multiplying by $n$ moves per sweep gives the physical time $n^3\log n/(8\pi^2)$.

\subsection{Architecture of the proof}\label{sec:overview}
The proof separates naturally into an upper-bound problem and a profile-identification problem.  The upper-bound problem is genuinely full-deck: it must control all $n!$ configurations and all irreducible components of the regular representation.  The profile problem is finer.  Once the full-deck remainder is known to be small at the $n^3$ window scale, one must identify the residual mode, prove its Gaussian limit, and then show that no other asymptotically detectable information survives.

The first simplification is geometric.  A Rudvalis move contains a deterministic cyclic displacement together with a single fair local randomization.  Observing the chain after complete rotations and passing to a rotating coordinate frame removes the deterministic drift.  One obtains the systematic adjacent-projection sweep
\[
 K=P_{n-1}\cdots P_0,
 \qquad P_i=\frac{I+U_{(i,i+1)}}2.
\]
This operator is not self-adjoint, so its powers cannot be studied by the usual reversible spectral theorem.  The role of the weighted interchange operator
\[
 L_Q=Q^{-1}A_0+\sum_{i=1}^{n-1}A_i
\]
is precisely to supply a reversible comparison object.  A local three-site calculation relates the Dirichlet energy of the endpoint $KF$ to the norm dissipated by the projection sweep.  From this one obtains a resolvent comparison for the singular values of $K$.  This is the mechanism that transports static heat-trace estimates for $L_Q$ back to the original directed shuffle without treating $K$ as if it were reversible.

The main full-deck difficulty is representation-theoretic multiplicity.  A bound for one card, or even for a fixed number of labelled cards, does not by itself control the regular representation.  We therefore group irreducible representations by the first-row deficit $r=n-\lambda_1$ of the corresponding partition.  For small and intermediate deficits, deletion-harmonic functions on ordered injective $r$-tuples contain enough copies of every relevant irreducible.  The comparison becomes effective only after smoothing the one-card kernel.  At the smoothed scale, zero extension from injective tuples to the full product incurs an explicit collision penalty of order $r^2p_*$, with $p_*$ the largest smoothed transition probability, instead of the exponentially large extension losses that arise at the original nearest-neighbor scale.  Young branching then pays the regular-representation multiplicities with the correct factorial factor.

Large deficits require a different idea.  A crude spectral gap is strong enough once the effective dimension has first been reduced.  We accomplish this with a short heat-trace entrance estimate obtained by cutting the path into blocks and restricting the regular representation to the corresponding Young subgroup.  This produces a sectorwise bound that depends on the deficit $r$, rather than paying the full dimension $n!$ at the final cutoff time.  Combining the harmonic-tuple estimate for moderate $r$ with the block entrance and rough gap for large $r$ closes the entire regular heat trace.

To work at the exact cutoff window, rather than at a fixed positive multiplicative slack, the weak-edge parameter $Q$ must tend slowly to infinity.  This removes the seam loss in the sweep-to-interchange comparison.  The required uniformity is supplied by a heat-kernel bound for the smoothed one-card chain that is independent of $Q$.  With a suitable choice of $Q$, smoothing scale, and deficit threshold, the full regular trace at
\[
 \frac{n^2}{2\pi^2}(\log n+2s)
\]
has the precise leading contribution coming from the two lowest one-card modes.  After transfer through the nonnormal resolvent estimate, this yields
\[
 \limsup_{n\to\infty}\chi_n^2(t_n(s))
 \le \exp(2\ee^{-2s})-1.
\]

The remaining task is to identify the random variable that exhausts this $L^2$ budget.  The one-card sweep has an exact eigenvalue determined by
\[
 \lambda^{n-2}(2\lambda-1)=1,
\]
and its associated complex Fourier eigenfunction lifts to a normalized statistic of the entire deck.  At time $t_n(s)$ this statistic retains a deterministic mean of order $\ee^{-s}$.  A one-sweep martingale decomposition reduces its fluctuations to bounded increments.  The key quantitative input is concentration of the conditional one-sweep covariance, proved by encoding the top--bottom interactions during a sweep in a random tree whose degree-square sum has expectation of order $n$.  The martingale central limit theorem then gives a Gaussian shift: under stationarity the relevant real statistic converges to $N(0,1)$, whereas from the ordered deck it converges to $N(\sqrt2\ee^{-s},1)$.

A one-dimensional Gaussian limit alone would give only a lower bound on the full-deck total-variation distance.  The sharp chi-square bound supplies the missing converse.  The Gaussian shift has likelihood ratio with second moment exactly $\exp(2\ee^{-2s})$, matching the full-deck limiting $L^2$ budget.  A bounded-truncation argument therefore forces the full likelihood ratio to become asymptotically measurable with respect to the single Gaussian statistic.  This saturation step rules out any additional surviving information and identifies the complete total-variation profile.

The only deep external comparison input is the Octopus Inequality~\cite{CLR}.  All other estimates needed to pass from that comparison to the Rudvalis profile are proved below.  Standard representation-theoretic facts about $\Sn$ are recalled explicitly and cited to Sagan~\cite{Sagan}.

\section{Definitions, notation, and background}\label{sec:prelim}

\subsection{Markov-chain distances and cutoff}
For probability measures $\mu,\nu$ on a finite set $\Omega$,
\[
 \norm{\mu-\nu}_{\TV}
 =\frac12\sum_{x\in\Omega}|\mu(x)-\nu(x)|
 =\max_{A\subseteq\Omega}|\mu(A)-\nu(A)|.
\]
If $\pi$ is strictly positive and $\mu\ll\pi$, the chi-square divergence is
\[
 \chi^2(\mu,\pi)=\sum_x\frac{(\mu(x)-\pi(x))^2}{\pi(x)}
 =\E_\pi\!\left[\left(\frac{d\mu}{d\pi}-1\right)^2\right].
\]
Cauchy--Schwarz gives
\[
 \norm{\mu-\pi}_{\TV}\le \frac12\sqrt{\chi^2(\mu,\pi)}.
\]
Both distances contract when the same Markov kernel is applied to both measures.

For $0<\varepsilon<1$, the total-variation mixing time is
\[
 t_{\mathrm{mix}}(\varepsilon)
 =\min\{t:\ d_n(t)\le\varepsilon\}.
\]
A sequence of chains exhibits total-variation cutoff at locations $t_n$ with windows $w_n=o(t_n)$ and profile $F$ if
\[
 d_n(\floor{t_n+s w_n})\longrightarrow F(s)
\]
for each fixed real $s$, where $F(-\infty)=1$ and $F(+\infty)=0$.  Standard background on mixing and cutoff may be found in~\cite{LPW,Lacoin2016}.

\subsection{The group action and the rotating-frame sweep}
Let $\Sn$ act on deck positions on the right.  For $g\in\Sn$, let $U_g$ denote the corresponding unitary operator on functions on any permutation representation:
\[
 (U_gf)(\sigma)=f(\sigma g).
\]
In particular, on $\C[\Sn]$ this is the right regular representation.  Every transposition $s$ satisfies $U_s^*=U_s=U_s^{-1}$.

Let $\rho$ be the cyclic positional shift that moves the top card to the bottom, and let $\tau$ interchange the bottom two positions.  One Rudvalis move is $\rho$ or $\tau\rho$, each with probability $1/2$; changing the convention for right versus left multiplication merely transposes the operator and does not affect any singular-value or mixing estimate below.  Thus, after fixing a convention, the one-step operator can be written as a deterministic rotation composed with one fair two-site average.

\begin{lemma}[One full rotation is an adjacent-projection sweep]\label{lem:sweepdecomp}
Let $s_i=(i,i+1)$ with indices understood modulo $n$, and set
\[
 P_i=\frac{I+U_{s_i}}2,\qquad A_i=I-P_i.
\]
After $n$ elementary Rudvalis moves, and after conjugating by the deterministic rotating frame, the Markov operator is
\[
 K=P_{n-1}P_{n-2}\cdots P_0.
\]
Consequently estimates for $K^m$ are estimates for the original shuffle after $nm$ elementary moves.  Reversing the orientation of the deck replaces $K$ by a conjugate or transpose, which has the same singular values.
\end{lemma}
\begin{proof}
Write one step as $U_\rho P$, where $P=(I+U_\tau)/2$ is the fair average at the bottom seam.  Using $U_\rho P U_\rho^{-1}$ to move the averaging operator through a rotation, expand
\[
 (U_\rho P)^n
 =U_\rho^n\,(U_\rho^{-(n-1)}PU_\rho^{n-1})\cdots(U_\rho^{-1}PU_\rho)P.
\]
Since $\rho^n=e$, the deterministic rotations cancel.  The conjugates of the bottom seam by successive powers of $\rho$ are exactly the $n$ adjacent transpositions around the cycle, in cyclic order.  Relabelling the starting edge gives the displayed product.
\end{proof}

\subsection{Linear-algebra conventions}
All finite-dimensional inner products are complex Hilbert-space inner products, linear in the second variable.  For self-adjoint operators $B,C$, the notation
\[
 B\preceq C
\]
means $\ip f{Bf}\le\ip f{Cf}$ for every vector $f$.  This is the Loewner or quadratic-form order.  If $0\preceq B\preceq C$, the min--max principle implies that the ordered eigenvalues of $B$ do not exceed the corresponding ordered eigenvalues of $C$.

For an operator $M$, $s_1(M)\ge s_2(M)\ge\cdots$ are its singular values, the square roots of the eigenvalues of $MM^*$.  The Hilbert--Schmidt norm is
\[
 \norm M_{\HS}^2=\Tr(MM^*)=\sum_j s_j(M)^2.
\]
We will use the Schatten--H\"older inequality in the simple form
\[
 \norm{M^m}_{\HS}^2\le \sum_j s_j(M)^{2m}.
\]
This is important because $K$ is not normal: one cannot replace singular-value arguments by powers of a Loewner inequality.

\subsection{The regular representation and partitions of $n$}\label{sec:repprimer}
A \emph{partition} $\lambda\vdash n$ is a nonincreasing sequence
$\lambda_1\ge\lambda_2\ge\cdots\ge0$ with sum $n$.  Its Young diagram has $\lambda_i$ boxes in row $i$; $\lambda'$ denotes the conjugate partition obtained by exchanging rows and columns.  The irreducible complex representations of $\Sn$ are indexed by partitions $\lambda\vdash n$.  We denote the irreducible module by $V_\lambda$ and its dimension by $d_\lambda$.

The regular representation decomposes as
\[
 \C[\Sn]\cong\bigoplus_{\lambda\vdash n} d_\lambda V_\lambda.
\]
Therefore, if an operator is built from the group action and hence respects this decomposition,
\[
 \Tr_{\C[\Sn]}T
 =\sum_{\lambda\vdash n}d_\lambda\Tr_{V_\lambda}T.
\]
The integer
\[
 r=n-\lambda_1
\]
will be called the \emph{first-row deficit}.  It measures how far $V_\lambda$ is from the trivial representation $(n)$ and is the natural parameter in the trace estimates below.

If $\nu\vdash r$, let $f^\nu$ denote the number of standard Young tableaux of shape $\nu$.  Equivalently, $f^\nu=\dim V_\nu$.  Young branching implies that the permutation representation on ordered injective $r$-tuples,
\[
 \C[\Omega_{n,r}]
 \cong \operatorname{Ind}_{\mathfrak S_{n-r}}^{\Sn}\ind,
\]
contains exactly $f^\nu$ copies of $V_{(n-r,\nu)}$.  We will also use
\[
 \sum_{\nu\vdash r}(f^\nu)^2=r!,
\]
which is simply the dimension identity for the regular representation of $\mathfrak S_r$.  These facts, the branching rule, and the content formula used later are standard; see~\cite{Sagan}.

For $1\le r\le n-1$, let
\[
 \Omega_{n,r}=\{(x_1,\ldots,x_r)\in[n]^r:\ x_a\ne x_b\text{ for }a\ne b\}
\]
be the ordered injective tuples.  For each coordinate $a$, the \emph{deletion marginal} $C_aF$ is obtained by summing $F$ over all allowed values of $x_a$, with the remaining coordinates fixed.  The deletion-harmonic space is
\[
 H_{n,r}=\bigcap_{a=1}^r\ker C_a.
\]
Each $C_a$ intertwines the action of $\Sn$, so $H_{n,r}$ is invariant.  The role of $H_{n,r}$ is to retain precisely enough copies of the deficit-$r$ irreducibles while forcing every coordinate to be mean zero after zero extension to $[n]^r$.

\subsection{A notation guide}
The symbols that recur most often are collected here.
\begin{center}
\begin{tabularx}{0.97\textwidth}{@{}p{0.20\textwidth}X@{}}
\toprule
Symbol & Meaning\\
\midrule
$R$ & one-elementary-move Rudvalis Markov operator.\\
$K$ & one full rotating-frame sweep, equal to $P_{n-1}\cdots P_0$ and corresponding to $n$ elementary moves.\\
$P_i,A_i$ & adjacent averaging projection $(I+U_{s_i})/2$ and its loss projection $I-P_i$.\\
$L_Q$ & self-adjoint weighted interchange operator with one weak cycle edge of weight $1/Q$.\\
$\ell_Q$ & one-card operator corresponding to $L_Q$.\\
$Z_{n,Q}(t)$ & regular-representation heat trace $\Tr_{\C[\Sn]}e^{-tL_Q}$.\\
$T_{n,r}(t)$ & part of that trace coming from partitions with first-row deficit $r$.\\
$H_{n,r}$ & deletion-harmonic functions on ordered injective $r$-tuples.\\
$\lambda_n$ & exact one-card eigenvalue of one elementary Rudvalis move near $e^{2\pi i/n}$.\\
$\Lambda_n=\lambda_n^n$ & corresponding eigenvalue after one sweep.\\
$Z(\sigma)$ & normalized complex deck statistic obtained by pairing label and position copies of the one-card eigenfunction.\\
\bottomrule
\end{tabularx}
\end{center}

\section{Symmetric comparison and transfer to the Rudvalis sweep}\label{sec:retained}
For cycle sites $\mathbb Z/n\mathbb Z$, put
\[
 s_i=(i,i+1),\quad P_i=(I+U_{s_i})/2,\quad A_i=I-P_i,\quad
 K=P_{n-1}\cdots P_0.
\]
For an integer $Q\geq1$, define
\[
 L_Q=Q^{-1}A_0+\sum_{i=1}^{n-1}A_i,\qquad
 c_Q=\frac{1+5/Q}{4},\qquad
 Z_{n,Q}(t)=\Tr_{\C[\Sn]}e^{-tL_Q}.
\]
Its one-card positive operator $\ell_Q$ has jump rates $1/2$ on ordinary edges and $1/(2Q)$ on the weak edge. Hence $P=I-\ell_Q$ is symmetric stochastic.

\subsection{Octopus smoothing and harmonic zero extension}
For a symmetric stochastic matrix $P$ put
\[
 \mathcal L(P)=\sum_{x<y}P(x,y)(I-U_{(xy)}).
\]
The following consequence of the octopus theorem~\cite{CLR} is also Proposition 3 of~\cite{HKS}. We use only its interchange comparison, not an exclusion cutoff theorem.

\begin{lemma}[Dynamical comparison]\label{lem:octopus}
For every $\tau\geq0$, in every unitary representation,
\[
 \mathcal L(e^{-\tau(I-P)})\preceq\tau\mathcal L(P).
\]
\end{lemma}
\begin{proof}
Write $D_{xy}=I-U_{(xy)}$ and $D_{xx}=0$. The octopus inequality at $o$, multiplied by the total off-diagonal incident weight, gives
\[
 (1-P(o,o))\sum_xP(o,x)D_{ox}
 \succeq \frac12\sum_{x,y\ne o}P(o,x)P(o,y)D_{xy}.
\]
Adding the terms with exactly one of $x,y$ equal to $o$ gives
\[
 \sum_xP(o,x)D_{ox}\succeq\frac12\sum_{x,y}P(o,x)P(o,y)D_{xy}.
\]
Sum in $o$ to get $2\mathcal L(P)\succeq\mathcal L(P^2)$. Apply this repeatedly to $(1-\tau/2^k)I+(\tau/2^k)P$ for $2^k\geq\tau$, and let $k$ increase. The octopus inequality in the regular representation implies it in every irreducible and hence every unitary direct sum.
\end{proof}

Let $\Omega_{n,r}$ be the ordered injective $r$-tuples of sites. Let $H_{n,r}$ be the common kernel of the $r$ deletion marginals, each marginal summing over the allowed values of the deleted coordinate. Use counting weight $n^{-r}$ on both $\Omega_{n,r}$ and the full product $[n]^r$. Zero extension $J_0$ is then isometric and maps $H_{n,r}$ into $(\ind^\perp)^{\otimes r}$. Deletion intertwines every site permutation, so $H_{n,r}$ is invariant for every interchange operator.

Set $p_\tau=e^{-\tau\ell_Q}$, $p_* =\max_{x,y}p_\tau(x,y)$, and
$\mathcal A_\tau=\sum_{a=1}^r(I-p_\tau)_a$ on the product. Let
$0=\gamma_0<\gamma_1\leq\cdots\leq\gamma_{n-1}$ be the one-card eigenvalues.

\begin{lemma}[Exact collision identity and trace bound]\label{lem:zero}
For every tuple function $F$,
\[
 \ip{J_0F}{\mathcal A_\tau J_0F}
 =\ip{F}{\mathcal L(p_\tau)F}-\mathcal E_{\rm swap}(F)
   +n^{-r}\sum_{x\in\Omega_{n,r}}V_\tau(x)|F(x)|^2,
\]
where $\mathcal E_{\rm swap}\geq0$ is the tracked-label exchange energy and
$V_\tau(x)=\sum_{a\ne b}p_\tau(x_a,x_b)\leq r(r-1)p_*$. Consequently,
\begin{equation}\label{eq:htrace}
 \Tr_{H_{n,r}}e^{-tL_Q}
 \leq \exp\!\left(\frac t\tau r(r-1)p_*\right)
 \left[\sum_{k=1}^{n-1}\exp\!\left\{-\frac t\tau(1-e^{-\tau\gamma_k})\right\}\right]^r.
\end{equation}
\end{lemma}
\begin{proof}
In the product Dirichlet form, a jump between two injective tuples is exactly a tracked-label-to-hole move. A jump from an injective tuple to a collision has terminal value zero; its two orientations cancel the form's factor $1/2$ and give the displayed killing term. Edges between collision tuples contribute zero. These observations give the identity and, by Lemma~\ref{lem:octopus}, the inequality
\[
 \ip{J_0F}{\mathcal A_\tau J_0F}
 \leq\tau\ip{F}{L_QF}+r(r-1)p_*\norm F^2.
\]
On $(\ind^\perp)^{\otimes r}$ the product spectrum is the multiset
$\xi=\sum_a(1-e^{-\tau\gamma_{k_a}})$ with every $k_a\geq1$.
If $\lambda_j$ are the ordered eigenvalues on $H_{n,r}$, min--max on the isometric image of the first $j$ eigenvectors gives
$\tau\lambda_j+r(r-1)p_*\geq\xi_j$. Exponentiate these scalar inequalities, sum, and enlarge the sum to the full product spectrum. This proves~\eqref{eq:htrace}.
\end{proof}

\subsection{Regular multiplicities and a rough sector gap}
For a partition $\lambda\vdash n$, write $d_\lambda=\dim V_\lambda$, and call
$r=n-\lambda_1$ its first-row deficit. Define
\[
 T_{n,r}(t)=\sum_{\lambda:n-\lambda_1=r}d_\lambda\Tr_{V_\lambda}e^{-tL_Q}.
\]
\begin{lemma}[Multiplicity and gap]\label{lem:multgap}
For $1\leq r\leq n-1$,
\begin{align}
 T_{n,r}(t)&\leq\binom nr\Tr_{H_{n,r}}e^{-tL_Q},\label{eq:mult}\\
 L_Q|_{V_\lambda}&\succeq \frac r{2n^2}I
 \quad\text{when }n-\lambda_1=r.\label{eq:gap}
\end{align}
\end{lemma}
\begin{proof}
Write $\lambda=(n-r,\nu)$. Young branching in
$\C[\Omega_{n,r}]=\operatorname{Ind}_{\mathfrak S_{n-r}}^{\Sn}\ind$
gives $f^\nu$ copies of $V_\lambda$. Every deletion target has first row at least $n-r+1$, so all these copies are harmonic. A standard tableau of shape $\lambda$ is determined injectively by the set of its $r$ entries below the first row and the standardized lower tableau. Thus
$d_\lambda\leq\binom nr f^\nu$. This proves~\eqref{eq:mult}.

For an ordinary $b$-site path set
\[
 L_b^{\rm path}=\frac12\sum_{i=1}^{b-1}(I-U_{(i,i+1)}),\qquad
 C_b=\sum_{i<j}(I-U_{(ij)}).
\]
Express $(i,j)$ as an adjacent word of length $2(j-i)-1$. Telescoping the unitary word and applying Cauchy--Schwarz gives a quadratic-form comparison. The total coefficient at edge $k$ is at most
\[
 4\sum_{i\leq k<j}(j-i)=2b\,k(b-k)\leq b^3/2.
\]
Since $\norm{(I-U_s)F}^2=2\ip F{(I-U_s)F}$, this yields
$C_b\preceq b^3L_b^{\rm path}$.
The central transposition scalar on $V_\lambda$ is
\[
 c_\lambda=\binom b2-\sum_i\binom{\lambda_i}2+
                      \sum_j\binom{\lambda'_j}2\geq br/2,
 \qquad r=b-\lambda_1.
\]
Indeed every row is at most $b-r$, so the negative row sum is at most $b(b-r-1)/2$. These standard branching and content identities are recorded, for example, in~\cite{Sagan}. Deleting the weak edge of $L_Q$ leaves an ordinary $n$-site path, proving~\eqref{eq:gap}.
\end{proof}

\begin{lemma}[Path trace]\label{lem:pathtrace}
For every $u\geq0$,
\[
 Z_b^{\rm path}(u):=\Tr_{\C[\mathfrak S_b]}e^{-uL_b^{\rm path}}
 \leq\exp\!\left(b^2e^{-u/(2b^2)}\right).
\]
\end{lemma}
\begin{proof}
The tableau bound and $\sum_{\nu\vdash r}(f^\nu)^2=r!$ give
$\sum_{\lambda:b-\lambda_1=r}d_\lambda^2\leq b^{2r}/r!$. Apply the gap in Lemma~\ref{lem:multgap}, include the trivial representation, and sum the exponential series. For $b=1$ the claimed bound also holds.
\end{proof}

\subsection{The sweep endpoint, resolvent, and Gamma transfer}
\begin{lemma}[Nonnormal-safe transfer]\label{lem:transfer}
For every integer $m\geq1$ and $0<\omega<1$,
\begin{align}
 K^*L_QK&\preceq c_Q(I-K^*K),\label{eq:endpoint}\\
 \chi_n^2(nm)&\leq\Tr(I+L_Q/c_Q)^{-m}-1\label{eq:resolvent}\\
 &\leq Z_{n,Q}((1-\omega)m/c_Q)-1+(n!-1)e^{-m\omega^2/2}.
 \label{eq:gamma}
\end{align}
\end{lemma}
\begin{proof}
Put $g_0=F$, $g_{i+1}=P_ig_i$, $d_i=A_ig_i$, and $G=KF$. Orthogonality at each individual projection gives
$D:=\norm F^2-\norm G^2=\sum_i\norm{d_i}^2$.
For neighboring transpositions $s,t$, let $C_{st}$ denote the orthogonal projector onto the sign representation of the three-site subgroup they generate, namely
\[
 C_{st}=\frac16\sum_{g\in\langle s,t\rangle}\operatorname{sgn}(g)U_g.
\]
A direct multiplication in the group algebra gives
\[
 A_sA_tA_s=\tfrac14 A_s+\tfrac34 C_{st}.
\]
This is verified directly in the six-element group algebra. If $P_sh=h$, the alternating part of $h$ vanishes, and the identity with $s,t$ reversed gives
$\norm{A_sA_th}^2=\norm{A_th}^2/4$.
For $1\leq i\leq n-2$, later disjoint projections commute with $A_i$, hence
\[
 A_iG=-P_{n-1}\cdots P_{i+2}A_id_{i+1}.
\]
The required preceding-average hypothesis holds for $h=g_{i+1}$. Thus
$\sum_{i=1}^{n-2}\norm{A_iG}^2\leq D/4$, while $A_{n-1}G=0$.
At the seam,
\[
 A_0G=-P_{n-2}\cdots P_2A_0d_1-A_0d_{n-1},\qquad
 \norm{A_0G}\leq\tfrac12\norm{d_1}+\norm{d_{n-1}},
\]
so $\norm{A_0G}^2\leq5D/4$. Weighting this edge by $1/Q$ proves~\eqref{eq:endpoint}.

The operator $(I+L_Q/c_Q)^{1/2}K$ is a contraction, implying
$KK^*\preceq(I+L_Q/c_Q)^{-1}$. Ordered eigenvalue comparison bounds each squared singular value. Schatten--H\"older then gives
\[
 \norm{K^m}_{\HS}^2\leq\sum_j s_j(K)^{2m}
 \leq\Tr(I+L_Q/c_Q)^{-m}.
\]
Convolution transitivity implies $1+\chi_n^2(nm)=\norm{K^m}_{\HS}^2$.
Finally, for $U_m$ Gamma with shape $m$ and rate one,
$(1+x)^{-m}=\E e^{-xU_m}$. Split the resulting expected trace at $(1-\omega)m$, use its monotonicity and $Z(0)=n!$, and apply
$\Pp(U_m<(1-\omega)m)\leq e^{-m\omega^2/2}$.
\end{proof}

\section{A sharp-window full-deck chi-square bound}\label{sec:sharpL2}
Two improvements remove the fixed-slack restriction without an unproved uniformity assertion.

\subsection{A heat-kernel bound uniform in the weak-edge weight}
\begin{lemma}[Uniform kernel and subdivision bounds]\label{lem:uniformkernel}
For every integer $Q\geq1$ and every $\tau>0$,
\begin{equation}\label{eq:uniformkernel}
 \max_{x,y}e^{-\tau\ell_Q}(x,y)\leq\frac1n+\sqrt{\frac{2\pi}{\tau}}.
\end{equation}
Also, with $N=n+Q-1$,
\begin{equation}\label{eq:subdivision}
 \gamma_j\geq2\sin^2\!\left(\frac{\pi\ceil{j/2}}N\right),
 \qquad 1\leq j\leq n-1.
\end{equation}
\end{lemma}
\begin{proof}
Delete the weak edge. On the remaining ordinary path choose a vertex where $|f|^2\leq\norm f_2^2/n$, with unnormalized counting norm. Along the unique path from that vertex to $x$, telescope $|f|^2$ and apply Cauchy--Schwarz:
\[
 |f(x)|^2\leq\frac{\norm f_2^2}{n}
       +2\norm f_2\left(\sum_{\rm path}|\nabla f|^2\right)^{1/2}
 \leq\frac{\norm f_2^2}{n}+2\sqrt2\norm f_2\sqrt{\ip f{\ell_Qf}}.
\]
The diagonal of the spectral projection onto eigenvalues at most $u$ is therefore at most $n^{-1}+2\sqrt{2u}$. Integration against $\tau e^{-\tau u}$ gives the diagonal bound~\eqref{eq:uniformkernel}. Positive semidefiniteness of the heat kernel bounds every off-diagonal entry by the geometric mean of the two diagonal entries. This argument never pays $Q$.

For~\eqref{eq:subdivision}, replace the weak edge by $Q$ ordinary edges. Linear interpolation on the inserted path preserves the one-card counting energy and increases counting norm. This injective map, applied to initial spectral subspaces including constants, gives the ordered eigenvalue comparison with the ordinary $N$-cycle.
\end{proof}

\subsection{A block entrance estimate in each deficit sector}
\begin{lemma}[Sectorwise entrance]\label{lem:sectorentrance}
For sufficiently large $n$, put
\[
 b=\floor{\frac n{\sqrt{12\log n}}},\qquad B=\ceil{n/b},\qquad u=n^2.
\]
Uniformly in $Q\geq1$ and $1\leq r\leq n-1$,
\begin{equation}\label{eq:sectorentrance}
 T_{n,r}(u)\leq\binom nr B^r e^{n^{-4}}.
\end{equation}
Consequently for $t\geq u$,
\begin{equation}\label{eq:tailsector}
 T_{n,r}(t)\leq
 e^{n^{-4}}\binom nr B^r\exp\!\left(-\frac{(t-u)r}{2n^2}\right).
\end{equation}
\end{lemma}
\begin{proof}
Partition the ordinary path into $B$ consecutive blocks of sizes $b_i\leq b$, and discard the interblock edges. On the ordered-tuple module, min--max gives
\[
 \Tr_{\Omega_{n,r}}e^{-uL_Q}\leq
 \Tr_{\Omega_{n,r}}e^{-uL_{\rm blocks}}.
\]
The block subgroup orbits are specified by assigning a block to each tracked label, so there are at most $B^r$ orbits. An orbit with occupancies $k_i$ is the product of the ordered $k_i$-tuple modules inside the blocks. Each such permutation module has irreducible multiplicities at most their dimensions, so its heat trace is at most the corresponding regular trace. Thus
\[
 \Tr_{\Omega_{n,r}}e^{-uL_{\rm blocks}}
 \leq B^r\prod_i Z_{b_i}^{\rm path}(u)
 \leq B^r\exp\!\left(\sum_i b_i^2e^{-n^2/(2b_i^2)}\right)
 \leq B^r e^{n^{-4}}.
\]
Here $\sum_i b_i^2\leq n^2$ and $n^2/(2b_i^2)\geq6\log n$. Combine with~\eqref{eq:mult} and the fact that the harmonic trace is at most the whole tuple trace. Finally split the exponential at $u$ in each deficit-$r$ eigenvalue and use~\eqref{eq:gap}.
\end{proof}

\subsection{Small, intermediate, and large deficits at the exact window}
\begin{proposition}[Sharp static upper profile]\label{prop:staticprofile}
Put $Q_n=\ceil{(\log n)^2}$. For any bounded sequence $s_n\to s$,
\[
 \limsup_{n\to\infty}\left[
 Z_{n,Q_n}\!\left(\frac{n^2}{2\pi^2}(\log n+2s_n)\right)-1\right]
 \leq\exp(2e^{-2s})-1.
\]
\end{proposition}
\begin{proof}
Write $L=\log n+2s_n$ and choose the explicit auxiliary parameters
\[
 a=n^{-1/250},\qquad \tau=\frac{a n^2}{2\pi^2},\qquad
 R_n=\floor{n^{99/100}},\qquad t=\frac{n^2L}{2\pi^2}.
\]
They are not claimed optimal. Lemma~\ref{lem:uniformkernel} gives
$p_*\leq n^{-1}(1+2\pi^{3/2}a^{-1/2})$. For $r\leq R_n$,
\begin{equation}\label{eq:collisionwindow}
 \frac1r\frac t\tau r(r-1)p_*
 \leq \frac L a\frac{R_n}{n}(1+2\pi^{3/2}a^{-1/2})
 =O(n^{-1/250}\log n)=o(1).
\end{equation}

The first two subdivision eigenvalues are at least
$(2\pi^2/n^2)(1-O(Q_n/n))$, while every eigenvalue with index at least three is at least
$(8\pi^2/n^2)(1-O(Q_n/n))$. Since $a\log n\to0$ and $(Q_n/n)\log n\to0$,
\begin{align*}
 z_n&:=\sum_{k=1}^{n-1}\exp\!\left[-\frac L a(1-e^{-\tau\gamma_k})\right]\\
 &\leq 2\exp(-L+o(1))+(n-3)\exp(-4L+o(1))
 =\frac{2e^{-2s}+o(1)}n.
\end{align*}
The $o(1)$ errors in the displayed exponents are absolute, not only relative to $L$. This follows from $1-e^{-x}=x+O(x^2)$ at the two lower spectral bounds; monotonicity handles every higher eigenvalue without expanding it individually.

Combining~\eqref{eq:htrace},~\eqref{eq:mult}, and~\eqref{eq:collisionwindow} gives a sequence $A_n\to2e^{-2s}$ such that, uniformly for $1\leq r\leq R_n$,
\[
 T_{n,r}(t)\leq\frac{A_n^r}{r!}.
\]
Therefore the sum in this entire range has limsup at most $\exp(2e^{-2s})-1$.

For $r>R_n$, use Lemma~\ref{lem:sectorentrance} and
$\binom nr\leq(en/r)^r$. Since $B\asymp\sqrt{\log n}$,
\[
 T_{n,r}(t)
 \leq e^{n^{-4}}
 \left[\frac{enB}{r}\exp\!\left(-\frac{\log n+2s_n}{4\pi^2}+\frac12\right)\right]^r.
\]
For these $r$, the logarithm of the bracket is at most
\[
 -\kappa\log n+\tfrac12\log\log n+O_s(1),\qquad
 \kappa=\frac1{4\pi^2}-\frac1{100}>0.
\]
For large $n$ it is at most $-\kappa\log n/2$. Summing the resulting bound over at most $n$ deficits gives a quantity tending to zero. Every nontrivial shape is covered.
\end{proof}

\begin{corollary}[Full-deck chi-square upper profile]\label{cor:L2}
For any integer sequence
\[
 m_n=\frac{n^2}{8\pi^2}(\log n+2s+o(1)),
\]
\begin{equation}\label{eq:chiupper}
 \limsup_n\chi_n^2(nm_n)\leq\exp(2e^{-2s})-1.
\end{equation}
In particular the full-deck densities have bounded $L^2(\pi)$ norm at every fixed window parameter.
\end{corollary}
\begin{proof}
Use Lemma~\ref{lem:transfer} with $Q_n=\ceil{(\log n)^2}$ and $\omega=n^{-1/4}$. Then
\[
 \frac{(1-\omega)m_n}{c_{Q_n}}
 =\frac{n^2}{2\pi^2}(\log n+2s+o(1)),
\]
because $(\log n)/Q_n\to0$ and $\omega\log n\to0$. The Gamma error is at most
$\exp[n\log n-\Theta(n^{3/2}\log n)]\to0$. Apply Proposition~\ref{prop:staticprofile}.
\end{proof}

\section{An exact eigenfunction and its elementary-move noise}\label{sec:eigen}
Here positions are numbered $1,\ldots,n$ as in the deck definition.

\begin{lemma}[Spectral root and eigenvector]\label{lem:root}
For all sufficiently large $n$ there is a root $\lambda_n=e^{z_n}$ of
\begin{equation}\label{eq:root}
 \lambda_n^{n-2}(2\lambda_n-1)=1
\end{equation}
with
\[
 z_n=\frac{2\pi i}{n}-\frac{4\pi^2}{n^3}+O(n^{-4}).
\]
Define $\Lambda_n=\lambda_n^n$, $\gamma_n=-\log|\Lambda_n|$, and choose
$\theta_n=\arg\Lambda_n$ close to zero. Then
\[
 \gamma_n=\frac{4\pi^2}{n^2}+O(n^{-3}),\qquad\theta_n=O(n^{-3}).
\]
There is an exact one-card eigenfunction $\phi_n$ with eigenvalue $\lambda_n$ such that
\begin{equation}\label{eq:phi}
 \sum_x\phi_n(x)=0,\quad \sum_x|\phi_n(x)|^2=n,\quad
 \max_x|\phi_n(x)|\leq2,
\end{equation}
and
\begin{equation}\label{eq:phismall}
 \left|\phi_n(n)-\phi_n(n-1)\right|=O(n^{-1}),\qquad
 \left|\sum_x\phi_n(x)^2\right|=O(1).
\end{equation}
\end{lemma}
\begin{proof}
On the analytic branch at zero,
\[
 (n-2)z+\log(2e^z-1)=nz-z^2+O(z^3).
\]
The equation with right side $2\pi i$ has a unique solution in a disk of radius $C/n^3$ about $2\pi i/n$: solve it by the contraction
$z\mapsto[2\pi i-\{\log(2e^z-1)-2z\}]/n$ on this disk, with $C$ sufficiently large. One substitution gives the stated expansion. In particular $|\lambda_n|<1$ eventually.

Before normalization set
\[
 \phi_n(x)=\lambda_n^{-(x-1)}\quad(1\leq x\leq n-1),\qquad\phi_n(n)=1.
\]
For positions $2,\ldots,n-1$ a card moves deterministically one place up. At positions $1$ and $n$ the one-card operator averages the values at $n-1,n$. Equation~\eqref{eq:root} therefore gives the eigenfunction relation at every site. The one-card matrix is doubly stochastic, so the eigenfunction has zero sum. Normalize its squared norm to $n$.

At the first $n-1$ sites it differs from $e^{-2\pi i(x-1)/n}$ by $O(n^{-2})$, before a normalization factor $1+O(n^{-2})$; at the last site the difference is $O(n^{-1})$. Thus its squared sum is $O(1)$ (in fact $O(n^{-1})$). Also
$\phi_n(n-1)=(2\lambda_n-1)\phi_n(n)$, giving the boundary difference bound. The sweep eigenvalue assertions follow by multiplying $z_n$ by $n$ and subtracting its $2\pi i$ phase.
\end{proof}

Suppress the subscript $n$ and put $a_j=\overline{\phi(j)}$. Define the complex deck statistic
\begin{equation}\label{eq:score}
 Z(\sigma)=\frac1{\sqrt n}\sum_{x=1}^n a_{\sigma(x)}\phi(x).
\end{equation}
It satisfies
\[
 RZ=\lambda Z,\qquad Z(\id)=\sqrt n.
\]
An elementary without-replacement calculation gives
\begin{equation}\label{eq:stationarymoments}
 \E_\pi Z=0,\qquad v_n:=\E_\pi|Z|^2=\frac n{n-1},\qquad
 q_n:=\E_\pi Z^2=\frac{|\sum_x\phi(x)^2|^2}{n(n-1)}=O(n^{-2}).
\end{equation}
For example the covariance of two distinct permuted $a$-values is minus their population second moment divided by $n-1$; expanding the two sums proves all formulas.

If $Z_t=Z(\sigma_t)$, let $\eta_{t+1}=Z_{t+1}-\lambda Z_t$. Conditional on the deck, the two possible outcomes differ by
\[
 \Delta(\sigma)=\frac{d_\phi}{\sqrt n}
       (a_{\sigma(1)}-a_{\sigma(n)}),\qquad
 d_\phi=\phi(n)-\phi(n-1).
\]
Hence the centered innovation is $\pm\Delta(\sigma)/2$, each sign with probability $1/2$, and
\begin{equation}\label{eq:micronoise}
 \E(|\eta_{t+1}|^2\mid\sigma_t)=\frac{|d_\phi|^2}{4n}|a_{\sigma_t(1)}-a_{\sigma_t(n)}|^2,
 \qquad |\eta_{t+1}|\leq Cn^{-3/2}.
\end{equation}
The corresponding conditional second moment, without conjugation, is obtained by squaring both complex differences.

\section{Concentration of the one-sweep conditional variance}\label{sec:bracket}
For a sweep starting at $\sigma$, put
\[
 \xi=Z(\sigma_n)-\Lambda Z(\sigma_0),\qquad
 b_n(\sigma)=\E(|\xi|^2\mid\sigma_0=\sigma),\qquad
 c_n(\sigma)=\E(\xi^2\mid\sigma_0=\sigma).
\]
We need concentration of these conditional quantities, not just their stationary means.

\begin{lemma}[Permutation edge-sum variance]\label{lem:edgesum}
Let $a_1,\ldots,a_n$ be bounded complex numbers, with a uniform bound independent of $n$. Let $G$ be a graph with at most $n$ indexed edges, allowing repetition, no loops, and degrees counting multiplicity. For deterministic complex weights $|w_e|\leq1$ and a uniform random permutation $\sigma$,
\[
 \E_\pi\left|\sum_e w_e\{h(a_{\sigma(i_e)},a_{\sigma(j_e)})-\bar h\}\right|^2
 \leq C\left(n+\sum_v\deg_G(v)^2\right),
\]
where either $h(u,v)=|u-v|^2$ or $h(u,v)=(u-v)^2$, and $\bar h$ is its common mean over two distinct population indices.
\end{lemma}
\begin{proof}
Each individual term is uniformly bounded. Pairs of edges sharing a vertex have bounded covariance, and their number is at most a constant times $\sum_v\deg(v)^2$. For two disjoint edges $e=\{i,j\}$ and $e'=\{k,\ell\}$, let $(X_1,X_2)$ and $(Y_1,Y_2)$ be two independent ordered pairs sampled uniformly without replacement from the population $\{a_1,\ldots,a_n\}$. Their joint law conditioned on the four sampled indices being distinct is exactly the law of $(a_{\sigma(i)},a_{\sigma(j)},a_{\sigma(k)},a_{\sigma(\ell)})$. The probability that the two pairs share an index is at most $4/(n-1)$. Since $h$ is uniformly bounded, conditioning changes any expectation of the product of the two edge observables by $O(1/n)$. Without the conditioning the two edge observables are independent, so the covariance for disjoint edges is $O(1/n)$. There are at most $n^2$ ordered pairs of disjoint edges. Expanding the squared absolute value of the weighted sum, and combining the disjoint-edge estimate with the shared-vertex count, proves the assertion. The same argument applies to complex weights and to the complex-valued observable $(u-v)^2$ by applying it to real and imaginary parts, or directly to complex covariance.
\end{proof}

\begin{lemma}[The graph exposed by one Rudvalis sweep]\label{lem:sweepgraph}
For a fixed sequence of fair insertion coins, view each top--bottom pair encountered during the first $n$ elementary moves as an edge between its two initial positions. This graph has $n$ edges, no loops, and
\[
 \E_{\rm coins}\sum_v\deg(v)^2\leq Cn.
\]
The graph depends on the coins and the initial positions, not on the actual card labels.
\end{lemma}
\begin{proof}
During the first $n-1$ moves, the top cards are successively the cards initially at positions $1,2,\ldots,n-1$. At move $j$ the bottom card is either the initial last card or one of the already processed top cards. Join the next new vertex to this bottom vertex. A bottom insertion makes the new vertex the bottom; a second-from-bottom insertion leaves the old bottom unchanged. Thus these $n-1$ edges form a tree with a bottom residence time dominated by a geometric random variable of success probability $1/2$.

For a nonroot vertex, its degree is one plus its number of later children. It acquires children only if its insertion made it bottom, an event of probability $1/2$; the ensuing child count is dominated by the same geometric variable, with first and second moments two and six. Its expected squared degree is therefore at most six. The root has the same upper bound. This gives at most $6n$ for the tree. The last move adds an edge between distinct existing vertices. Its increase of the degree-square sum is at most $4n$, which suffices. All these statements concern positions and coins alone.
\end{proof}

\begin{proposition}[Conditional covariance concentration]\label{prop:bracket}
Uniformly over decks,
\begin{equation}\label{eq:uniformbc}
 |\xi|\leq Cn^{-1/2},\qquad 0\leq b_n\leq Cn^{-2},\qquad |c_n|\leq Cn^{-2}.
\end{equation}
Their stationary means are
\begin{equation}\label{eq:bcmeans}
 \bar b_n=(1-|\Lambda|^2)v_n,\qquad
 \bar c_n=(1-\Lambda^2)q_n.
\end{equation}
Moreover,
\begin{equation}\label{eq:bcL2}
 \norm{b_n-\bar b_n}_{L^2(\pi)}+
 \norm{c_n-\bar c_n}_{L^2(\pi)}=O(n^{-5/2}).
\end{equation}
\end{proposition}
\begin{proof}
Iteration of the exact one-step eigenfunction relation gives
\[
 \xi=\sum_{j=0}^{n-1}\lambda^{n-1-j}\eta_{j+1}.
\]
The uniform innovation bound proves the first inequality in~\eqref{eq:uniformbc}. Conditional orthogonality of distinct martingale innovations gives
\begin{align}
 b_n(\sigma)&=\frac{|d_\phi|^2}{4n}\sum_{j=0}^{n-1}|\lambda|^{2(n-1-j)}
 \E\bigl[|a_{\sigma_j(1)}-a_{\sigma_j(n)}|^2\mid\sigma_0=\sigma\bigr],\label{eq:bformula}\\
 c_n(\sigma)&=\frac{d_\phi^2}{4n}\sum_{j=0}^{n-1}\lambda^{2(n-1-j)}
 \E\bigl[(a_{\sigma_j(1)}-a_{\sigma_j(n)})^2\mid\sigma_0=\sigma\bigr].\label{eq:cformula}
\end{align}
The remaining uniform bounds follow because the prefactors are $O(n^{-3})$ and there are $n$ bounded summands. At stationarity, expansion of
$Z_n=\Lambda Z_0+\xi$ and conditional centering of $\xi$ give~\eqref{eq:bcmeans}.

For~\eqref{eq:bcL2}, condition on all sweep coins and apply Lemma~\ref{lem:edgesum} to the graph in Lemma~\ref{lem:sweepgraph}. The weights in~\eqref{eq:bformula} or~\eqref{eq:cformula} are deterministic and have modulus at most one. The conditional mean over the uniform initial permutation is independent of the graph, because every edge joins two distinct positions. Average the variance bound over coins to obtain $O(n)$ for the unscaled sum. Conditional expectation over coins is an $L^2$ contraction. Multiplying by the $O(n^{-3})$ prefactor yields a standard deviation $O(n^{-5/2})$, as asserted.
\end{proof}

\section{A Gaussian shift at the cutoff window}\label{sec:clt}
Take any integer sequence
\[
 m_n=\frac{\log n/2+s+o(1)}{\gamma_n},\qquad s\in\R,
\]
and define the real statistic at this observation time by
\begin{equation}\label{eq:W}
 W_{n,m}(\sigma)=\sqrt2\,\Re\bigl(e^{-im\theta_n}Z(\sigma)\bigr).
\end{equation}
The phase is deterministic and is part of the observable, not extra randomization of the shuffle.

\begin{proposition}[Both marginal Gaussian limits]\label{prop:clt}
Under $\pi$ and under the shuffle law after $nm_n$ elementary moves, respectively,
\[
 W_{n,m_n}\ \Longrightarrow\ N(0,1),\qquad
 W_{n,m_n}\ \Longrightarrow\ N(\sqrt2e^{-s},1).
\]
\end{proposition}
\begin{proof}
Observe the chain every $n$ moves and write $\sigma^{[k]}=\sigma_{kn}$. Let
$\xi_k=Z(\sigma^{[k]})-\Lambda Z(\sigma^{[k-1]})$.
Then
\[
 W_{n,m}(\sigma^{[m]})
 =\sqrt2\,\Re(e^{-im\theta}\Lambda^m Z(\sigma^{[0]}))+
   \sum_{k=1}^m X_{n,k},
\]
where
$X_{n,k}=\sqrt2\,\Re(e^{-im\theta}\Lambda^{m-k}\xi_k)$
is a real martingale difference. Its maximum absolute value is at most $Cn^{-1/2}$.
The predictable variance sum is
\begin{align}\label{eq:V}
 V_n={}&\sum_{k=1}^m|\Lambda|^{2(m-k)}b_n(\sigma^{[k-1]})\\
 &+\Re\sum_{k=1}^m e^{-2im\theta}\Lambda^{2(m-k)}c_n(\sigma^{[k-1]}).
 \nonumber
\end{align}
It is nonnegative and uniformly bounded, since~\eqref{eq:uniformbc} and
$1-|\Lambda|^2\asymp n^{-2}$ bound both absolute sums.

\paragraph{Replacing conditional variances near the observation time.}
Fix $A>0$ and set $M=\ceil{A/\gamma_n}$. For the late indices $k>m-M$, the density of $\sigma^{[k-1]}$ has bounded $L^2(\pi)$ norm depending only on $s,A$. Indeed its earliest time is
$(\log n/2+s-A+o(1))/\gamma_n$, to which Corollary~\ref{cor:L2} applies; $L^2$ distance subsequently contracts under the Markov operator. Cauchy--Schwarz and~\eqref{eq:bcL2} give
\[
 \E|b_n(\sigma^{[k-1]})-\bar b_n|+
 \E|c_n(\sigma^{[k-1]})-\bar c_n|
 \leq C_{s,A}n^{-5/2}.
\]
Summing against the geometric weights in~\eqref{eq:V} gives total expected error at most $C_{s,A}n^{-1/2}\to0$.
For the early indices $k\leq m-M$, the uniform bounds~\eqref{eq:uniformbc} give a total absolute contribution at most $Ce^{-2A}$, and the same bound holds for the corresponding stationary means. Thus first taking $n\to\infty$ and then $A\to\infty$ justifies replacement by stationary means in $L^1$.

For a stationary initial state, the same replacement holds directly for all indices, without the late-time density argument. The replacement of the $b$ sum equals
\[
 \bar b_n\sum_{k=1}^m|\Lambda|^{2(m-k)}
 =v_n(1-|\Lambda|^{2m})\longrightarrow1.
\]
The absolute value of the replacement $c$ sum is at most
$|\bar c_n|/(1-|\Lambda|^2)=o(1)$, by~\eqref{eq:bcmeans},~\eqref{eq:stationarymoments}, and
$|1-\Lambda^2|=O(n^{-2})$. Consequently $V_n\to1$ in $L^1$ under either initial law.

\paragraph{An elementary bounded-increment martingale CLT.}
For completeness, if real martingale differences $X_k$ have maximum size $\delta_n\to0$, predictable variance sum $V_n\leq C$, and $V_n\to1$ in probability, then their sum tends to $N(0,1)$. To see this, fix real $u$, put $S_j=\sum_{k\leq j}X_k$, $V_j=\sum_{k\leq j}\E(X_k^2\mid\mathcal F_{k-1})$, and expand conditionally:
\[
 \E(e^{iuX_k}\mid\mathcal F_{k-1})
 =1-\frac{u^2}2 v_k+O_u(\delta_n v_k),
\]
because $\E(|X_k|^3\mid\mathcal F_{k-1})\leq\delta_n v_k$. Multiplying by $e^{u^2v_k/2}=1+u^2v_k/2+O_u(v_k^2)$ gives a one-step compensated error $O_u(\delta_n v_k+v_k^2)$. Since $v_k\leq\delta_n^2$, $\sum_k v_k\leq C$, and therefore $\sum_kv_k^2\leq C\delta_n^2$, telescoping expectations yields
$\E\exp(iuS_m+u^2V_m/2)=1+o(1)$. The compensated exponentials are uniformly bounded in modulus by $e^{u^2C/2}$. Boundedness and convergence of $V_m$ give the characteristic function limit $e^{-u^2/2}$.
Applying this result to~\eqref{eq:V} proves the centered Gaussian limit.

From the identity start, the initial term is
$\sqrt{2n}|\Lambda|^m\to\sqrt2e^{-s}$. From stationarity it tends to zero in $L^2$, since $\E_\pi|Z|^2$ is bounded and $|\Lambda|^{2m}=O(1/n)$. This proves both assertions.
\end{proof}

\begin{remark}[No circular mixing premise]
The late-time density bound comes solely from the independently proved full-deck chi-square upper bound. It does not assume the profile, a Gaussian approximation, or a coupling to stationarity. The constant $C_{s,A}$ may be very large; the limits are taken with $A$ fixed first. This is enough for the $n^{-1/2}$ error to vanish.
\end{remark}

\section{From the scalar limit to the full-deck profile}\label{sec:squeeze}
\begin{lemma}[Gaussian likelihood-ratio saturation]\label{lem:saturation}
Let $Q_n\ll P_n$ with likelihood ratio $L_n$, and let $Y_n$ be a real statistic. Suppose, for fixed $a\in\R$,
\[
 Y_n\Longrightarrow N(0,1)\quad(P_n),\qquad
 Y_n\Longrightarrow N(a,1)\quad(Q_n),\qquad
 \limsup_n\E_{P_n}L_n^2\leq e^{a^2}.
\]
Then
\[
 \norm{Q_n-P_n}_{\TV}\longrightarrow2\Phi(|a|/2)-1,
 \qquad \E_{P_n}L_n^2\longrightarrow e^{a^2}.
\]
Moreover $L_n$ under $P_n$ converges in distribution to
$\exp(aG-a^2/2)$, where $G\sim N(0,1)$.
\end{lemma}
\begin{proof}
The Gaussian likelihood ratio is $\ell(y)=e^{ay-a^2/2}$. For $M>0$ let
$g_M(y)=\min\{M,\ell(y)\}$, a bounded continuous function. The two weak limits and the squared-norm bound imply
\begin{align*}
 \limsup_n\E_{P_n}(L_n-g_M(Y_n))^2
 &\leq e^{a^2}-2\E_{N(a,1)}g_M+\E_{N(0,1)}g_M^2\\
 &=\E_{N(0,1)}(\ell-g_M)^2\longrightarrow0
 \quad(M\to\infty).
\end{align*}
Thus $L_n$ is approximated in the iterated $L^2$ sense by these bounded functions of the statistic. The inequality
\[
 \left|\frac12\E_{P_n}|L_n-1|-\frac12\E_{P_n}|g_M(Y_n)-1|\right|
 \leq\frac12\norm{L_n-g_M(Y_n)}_{L^2(P_n)}
\]
gives the TV limit, because the Gaussian TV distance is $2\Phi(|a|/2)-1$.
The same approximations give convergence in distribution of $L_n$ and the lower bound
$\liminf\norm{L_n}_2\geq\norm\ell_2$, matching the hypothesis.

This proof uses bounded truncations. It does not silently assume convergence of exponential moments of $Y_n$, nor a direct untruncated $L^2$ approximation by $\ell(Y_n)$.
\end{proof}

\begin{proof}[Proof of Theorem~\ref{thm:mainprofile}]
First take sweep times with
$m=(\log n/2+s+o(1))/\gamma_n$. Use $P_n=\pi$, $Q_n=\mu_{nm}$,
$Y_n=W_{n,m}$ and $a=\sqrt2e^{-s}$. Corollary~\ref{cor:L2} gives
$\limsup\E_\pi L_n^2\leq e^{2e^{-2s}}=e^{a^2}$, and Proposition~\ref{prop:clt} gives the two marginal Gaussian limits. Lemma~\ref{lem:saturation} proves both profiles at sweep times.

The root expansion implies
\[
 \frac{n}{\gamma_n}\left(\frac{\log n}2+s\right)
 =\frac{n^3}{8\pi^2}(\log n+2s)+O(n^2\log n),
\]
and the error is $o(n^3)$. More precisely the desired integer elementary time lies between neighboring sweep times whose normalized parameters tend to the same $s$. Both TV and chi-square contract under further elementary moves. Squeeze between these sweep times; the proofs above allow $s_n\to s$. This establishes~\eqref{eq:main-tv}--\eqref{eq:main-chi}.

For the likelihood-ratio statement below, we also record the scalar limit at nonsweep times rather than relying only on monotonicity. Write $t=nm+j$, $0\leq j<n$, choose $\alpha_n=\arg\lambda_n$ with $n\alpha_n=2\pi+\theta_n$, and set
\[
 W_{n,t}^{\rm elem}(\sigma)=\sqrt2\,\Re(e^{-it\alpha_n}Z(\sigma)).
\]
Iterating the elementary innovations over the final $j$ moves gives, on every trajectory,
\[
 \left|e^{-it\alpha_n}Z(\sigma_t)
       -|\lambda_n|^j e^{-inm\alpha_n}Z(\sigma_{nm})\right|
 \leq Cn^{-1/2}.
\]
Also $1-|\lambda_n|^j=O(n^{-2})$ and $|Z(\sigma)|\leq C\sqrt n$ uniformly. Therefore $W_{n,t}^{\rm elem}(\sigma_t)$ and $W_{n,m}(\sigma_{nm})$ differ by $O(n^{-1/2})$ uniformly. Under either identity or stationary initialization the same two Gaussian limits follow. Lemma~\ref{lem:saturation}, with the chi-square upper bound just obtained by squeezing, consequently applies directly at every stated elementary time as well.
\end{proof}

\begin{corollary}[Log-likelihood and relative entropy]
At the times of Theorem~\ref{thm:mainprofile}, the likelihood ratio under $\pi$ converges in law to
\[
 \exp\!\left(\sqrt2e^{-s}G-e^{-2s}\right).
\]
Its logarithm therefore converges in law to $N(-e^{-2s},2e^{-2s})$. Also
\[
 D(\mu_{t_n(s)}\Vert\pi)\longrightarrow e^{-2s}.
\]
This is a window-scale entropy profile, not a uniform nonlinear entropy contraction theorem.
\end{corollary}
\begin{proof}
Use Lemma~\ref{lem:saturation}. The limiting likelihood ratio is strictly positive, so the logarithm is a permissible continuous mapping, with any zero likelihoods treated as $-\infty$. For relative entropy, $x\log x$ is bounded below and grows strictly slower than $x^2$ at infinity. The uniformly bounded second moments give uniform integrability of $L_n\log L_n$, and its limiting expectation is $a^2/2=e^{-2s}$.
\end{proof}

\section{Consequences, comparison with earlier work, and reproducibility}\label{sec:discussion}

\subsection{What the theorem adds to the classical picture}
Hildebrand's upper bound and Wilson's lower bound established the order $\Theta(n^3\log n)$, with Wilson identifying the lower constant $1/(8\pi^2)$~\cite{Hildebrand1990,Wilson2003}.  Theorem~\ref{thm:mainprofile} supplies the matching upper constant, proves cutoff, identifies the $n^3$ window, and determines the entire total-variation and chi-square profiles.  In particular, the leading-order lower bound no longer has to be imported from Wilson: it follows from the negative-$s$ side of the profile itself.  Wilson's theorem remains the historical source of the sharp coefficient on the lower-bound side.

The mechanism also explains the logarithmic separation between one-card and full-deck relaxation.  The slowly decaying Fourier mode has sweep contraction
$|\Lambda_n|=\exp\{-4\pi^2/n^2+O(n^{-3})\}$.  Thus one copy relaxes on the order-$n^2$ sweep scale, equivalently $n^3$ elementary moves.  The ordered deck begins with a coherent sum of $n$ copies, of size $\sqrt n$ after normalization.  Reducing this signal to order one requires an additional factor $\frac12\log n$ in the exponent, producing the $n^3\log n$ mixing scale.

The limiting likelihood ratio under stationarity is lognormal:
\[
 L_s=\exp\!\left(\sqrt2\ee^{-s}G-\ee^{-2s}\right),\qquad G\sim N(0,1).
\]
Its second moment is $\exp(2\ee^{-2s})$, which matches the limiting chi-square divergence plus one.  Its total-variation distance from one is the Gaussian-shift distance in \eqref{eq:main-tv}.  This equality of the full-deck $L^2$ budget with the one-dimensional Gaussian likelihood is the reason a single statistic determines the entire profile here.

\subsection{Relation to exclusion and interchange methods}
The smoothing step is closely related in spirit to the recent use of the Octopus Inequality in exclusion-process comparison~\cite{HKS}.  There is, however, an important distinction.  The Rudvalis sweep is nonnormal and is not itself an exclusion or interchange process.  The reversible interchange operator $L_Q$ is an auxiliary comparison object.  The proof therefore has two additional layers: a nonnormal endpoint/resolvent transfer from $K$ to $L_Q$, and a full regular-representation trace estimate that pays all labelled multiplicities.  These layers are what allow a reversible spectral estimate to control the original directed shuffle.

The representation-theoretic organization is also different from a fixed-density exclusion problem.  The regular representation contains every partition $\lambda\vdash n$ with multiplicity $d_\lambda$.  The deletion-harmonic tuple spaces cover all small-to-intermediate first-row deficits with the correct factorial payment, while the sectorwise entrance estimate handles the remaining large deficits without ever paying the raw dimension $n!$ at the cutoff window.

\subsection{Computational checks}
The proof is analytic and does not rely on numerical certification.  For reproducibility, the ancillary script \texttt{checks/verify\_rudvalis\_profile.py} checks the finite-dimensional identities used in the argument on small symmetric groups and tuple modules.  It verifies, among other items, the collision-energy identity, the octopus-derived matrix comparison, the endpoint and resolvent inequalities, the exact one-card eigenfunction, and the conditional covariance formulas.  These checks are useful for normalization and sign errors but are not substitutes for the all-$n$ proofs above.

\section{Acknowledgments}
The author thanks Persi Diaconis for suggesting the problem, amongst many other interesting card shuffling problems. Gpt-6 astra pro was used extensively for derivations and polishing.

\end{document}